# An Optimization Approach to Improve Equitable Access to Local Parks

Anisa Young [a], Emily L. Tucker [a,*], Mariela Fernandez [b], David White [c], Robert Brookover [c], Brandon Harris [d]

[a] Department of Industrial Engineering, Clemson University, 271 Freeman Hall, Clemson, SC, USA

[b] Department of Recreation, Sport and Tourism, University of Illinois Urbana-Champaign, 1206 S 4th St., 226 Huff Hall, Champaign, IL, USA

[c] Department of Parks, Recreation, and Tourism Management, Clemson University, 142 Jersey Lane, Clemson, SC, USA

[d] Department of Teaching, Learning, and Sociocultural Studies, University of Arizona, 1430 E. Second Street, Tucson, AZ, USA

**Abstract**

Local parks are public resources that promote human and environmental welfare. Unfortunately, park inequities are commonplace as historically marginalized groups may have insufficient access. Platforms exist to identify the geographical areas that would benefit from future park improvements. However, these platforms do not optimize decisions nor include key features, such as budget and infrastructure, that are relevant to park location decisions. To support recreational and government agencies in addressing inequities in the distribution and quality of parks, we propose a mixed-integer program that minimizes insufficient access, defined as weighted deviations across multiple categories (distance, capacity, and environmental features). We consider an equity-focused min-max objective and an overall objective to minimize total weighted deviations. We apply the model to a case study of Asheville, North Carolina. We conduct extensive data collection to parameterize the model. In policy analyses, we consider the effects of available budget, planning horizons, strategic demographic priorities, and thresholds of access. The model reflects user-defined criteria and goals, and the results suggest that the framework may be generalizable to other cities. This study serves as a step in the development and incorporation of mathematical modeling to achieve social goals within the recreational setting.

Email addresses: anisay@clemson.edu (Young), etucke3@clemson.edu (Tucker), mfrnndz2@illinois.edu (Fernandez), whitedl@clemson.edu (White), bob@clemson.edu (Brookover), brandonharris@arizona.edu (Harris)

* Corresponding author. Current address: 271 Freeman Hall, Clemson, SC 29634

# 1. Introduction

## 1.1. Context and motivation

Local parks are foundational to healthy communities. They support local residents' physical, mental, and social health (Bedimo-Rung et al., 2005; Grahn and Stigsdotter, 2010; Kaczynski and Henderson, 2007; Miles et al., 2012; Rigolon, 2016). Their publicly available, open spaces promote physical activity, which may lead to decreases in ailments such as heart disease, cancer, and obesity (Bedimo-Rung et al., 2005; Cleven et al., 2020; Müller-Riemenschneider et al., 2020). Parks also improve the quality of human life by decreasing air and noise pollution, assisting with water runoff, and regulating temperature (Rigolon, 2016).

Unfortunately, there are persistent racial and ethnic disparities in access to local parks (Rigolon, 2016). Racially marginalized and impoverished communities are more likely to experience excessive distances to parks, overcrowding, and have parks with poor maintenance, crime, and insufficient amenities (Rigolon, 2016; Stodolska and Shinew, 2010). Park disparities stem from complex historical and political factors that work to undermine marginalized groups; marginalization is defined by non-white racial-ethnicity, low-income economic status, age dependency, and physical or mental disability (Rigolon, 2016; Rigolon et al., 2019). Historic discrimination and negative attitudes of distrust couple with other sentiments created by "neighborhood stigma" which can lead to separation and fear between racial groups (Harris et al., 2021). Gentrification, a social process by which marginalized groups must abandon their homes as a result of increased rental prices (Rigolon and Németh, 2018), can be exacerbated due to park development and may force marginalized groups to abide within underdeveloped and poorly maintained areas, many of which do not incorporate space for parks and greenspaces.

According to the environmental justice principles, there should be an equitable distribution of and fair access to natural environments and sufficient access to the aforementioned benefits of parks and greenspaces. This field also prioritizes the integration of marginalized groups into park planning processes (Rigolon et al., 2019). Finally, the enjoyment and use of park spaces should be inclusive and welcoming for all, regardless of race, gender, age, economic status, and ability-level.

When determining the effectiveness of parks to provide an appropriate social and recreational space for the community, there are several measures that can be considered. Quantitative measures include park accessibility, distribution, quantity, and capacity (Boulton et al., 2018; Hunter and Luck, 2015; Rigolon, 2017). Qualitative components relate to park maintenance, visual appeal, safety of the park, community preferences, and opportunities for mental and physical human welfare (Hunter and Luck, 2015; McCormack et al., 2010). Individual access and quality metrics may also be translated into an aggregate measure of access. One non-profit organization, the Trust for Public Land uses a scoring system that incorporates factors of population density, income, racial composition, health, heat, and air pollution (The Trust for Public Land, 2022a). The total score is the sum of each factor, weighted equally. The tool provides a map of the top general locations to place future parks considering underserved areas. However, it does not optimize the park plan, and it does not include a budget or community-defined scoring weights.

Taken together, much of the existing research on park disparities is evaluative, rather than *proactive*. The focus has been on understanding existing and historical inequities (Hughey et al., 2016; Rigolon, 2016, 2017; Stodolska and Shinew, 2010; Taylor et al., 2007). It is within this space of proactive planning that mathematical optimization can provide value to park practitioners. We define the problem of park access in a discrete facility location context (Daskin, 2013); a planner seeks to decide which eligible park locations (land parcels) to purchase to provide park access for local residents. We situate the problem in this context because location theory is commonly used by urban planners to inform development decisions, and the resulting Geographic Information System (GIS) maps are easily understandable to other municipal partners and community residents (Bara et al., 2022; Ibes, 2015; Talen, 1998).

Facility location models have been used to recommend strategic locations for several types of urban services (Farahani et al., 2019), companies (Manzini and Gebennini, 2008), healthcare sites (Pourrezaie-Khaligh et al., 2022), pre-positioning supplies for disaster relief (Alem et al., 2021; Shehadeh and Tucker, 2022), and conservation sites (Billionnet, 2013; Messer et al., 2016). Similar discrete location models are relevant to park planning, but there are modeling challenges that are unique to the outdoor recreational setting. Environmental factors, including tree cover and heat, affect park quality (similar to conservation models), though desirable park locations are those closer to resident rather than far (in contrast to conservation models and similar to for-profit models). Equity is key, and balancing needs across resident demographics is necessary (similar to public sector models), though defining capacity metrics by demographic group introduces nonlinearities. Finally, extensive geospatial and demographic data is available, and several datasets, with heterogeneous parameters and characteristics, need to be integrated. We seek to develop a model to address these gaps and provide a foundation for further research in the recreational setting.

**1.2. Related Literature**

The research on discrete facility location is wide (reviews available in Farahani et al., 2010; Melo et al., 2009; Owen & Daskin, 1998), and common objective functions are coverage, center, and median. Coverage models seek to locate facilities such that demand nodes are within a threshold distance (Berman et al., 2010). These have been applied to contexts ranging from ambulance services (van den Berg et al., 2016) to hand-sanitizing stations (O'Brien et al., 2021). P-center problems focus on the worst-case and minimize the distance to the customer that is furthest away (e.g., refueling stations, Lin and Lin, 2018). The p-median problem minimizes the total weighted distance to all customers (Daskin and Maass, 2015), e.g., minimizing demand-weighted distance to healthcare facilities (Jia et al., 2014). Trade-offs can be also considered between median objectives and equity-focused goals, e.g., facility workload (Daskin and Tucker, 2018).

Equity-focused models in this space have generally considered the differences in these distances across customers (Barbati and Piccolo, 2016). Equity is then measured as a function of the metric, e.g., center, range, mean absolute deviation, variance, and maximum deviation (Barbati and Piccolo, 2016; Marsh and Schilling, 1994). It has often been studied in the context of particular application areas. Within healthcare, for example, Shehadeh and Snyder (2022) provide a review of models of equity in stochastic location problems. Zhang et al. (2016) present a case study of locating healthcare facilities in Hong Kong; they consider equity, overall accessibility, coverage, and cost. Pourrezaie-Khaligh et al. (2022) present a healthcare facility location model that defines equity as a function of accessibility deviations, including capacity, demand, and travel time. Certain populations may also have greater need for services. Alem et al. (2021) apply the Social Vulnerability Index to consider equity in the preparation and response to disasters in Brazil. Similarly, in ambulance dispatching, priority is given according to the severity of health risk, and demand types are modeled with different priority indices (Enayati et al., 2019).

In the context of recreation, there are a handful of optimization models have been used to recommend park locations. Li (2014) presents maximum coverage and capacitated location-allocation models to meet greenspace standards of distances and acres per person. She provides a case study of Luohu District, Shenzhen, China. Sefair et al., (2012) develop a multi-objective location model for parks in Bogotá, Colombia. The model considers coverage, cost, connectivity, and residents who benefit and is solved using lexicographic methods. Leboeuf et al., (2023) present a two-stage approach to park planning in Montreal. In the first-stage, budgets are fairly allocated to specific neighborhoods, and in the second, the park locations are optimized and designed to maximize a function of expected park visits. Continuous location models optimize parks with multiple objectives (Yuan et al., 2011) and may be solved with meta-heuristics such as Genetic Algorithms (Neema and Ohgai, 2010). In addition to informing design, optimization may also be used to evaluate greenspace access; Niedermann et al., (2018) use a network flow model to find the minimum average resident distance to greenspace subject to capacity limits. These papers provide a foundation for optimization in greenspace planning. Yet, existing models do not consider equity in park access across individual subpopulations. This nuance is critical to address underlying park disparities and requires further modeling of access.

Related literature has considered optimization for land-use in urban settings more broadly, including urban development (Gabriel et al., 2006; Türk & Zwick, 2019), use allocation (Ligmann-zielinska et al., 2005; Liu et al., 2012; Shaygan et al., 2014), and mitigation of sprawl (Kumar et al., 2016). These styles of models tend to focus on city composition, land suitability, and compactness. In contrast, park planning focuses on providing access to residents, a person-centered rather than a land-centered objective. Farahani et al. (2019) review discrete location models for urban services and categorize existing research into six groups: waste management, large-scale disasters, small-scale emergency, general services and infrastructure, non-emergency health care, and transportation infrastructure. Park planning would fit within the general services and infrastructure category, though requires capacity, assignment, and budget considerations that have not been addressed together.

Optimization has been used in other outdoor-focused applications as well, including environmental conservation and invasive species management. A review by Billionnet (2013) provides an overview of several models to support biodiversity. The purpose of conservation sites is to preserve unique species of vegetation and animals (Boulton et al., 2018), and complexities include compactness and connectivity (Billionnet, 2013; Hamaide et al., 2022). These allow for animal and plant species to travel between connected sites. Öhman and Lämås (2005) study compact reserves for forest planning. St. John et al., (2018) develop methods to consider the geometry of the corridors through path planning and network optimization, and Álvarez-Miranda et al., (2021) consider both connectivity and buffer zones. Messer et al. (2016) introduce a multiple-knapsack structure to maximize the overall benefit of conservation program outcomes given capacity and budget constraints. To mitigate the appearance and spread of pests, Büyüktahtakın and Haight (2017, 2018) consider both spatial and temporal components in resource allocation; the goal is to minimize destruction caused by invasive species. Epanchin-Niell and Wilen (2012) develop optimal control policies to reduce long-term costs.

### 1.3. Contribution

In this study, we seek to address the problem of equitable, strategic park planning by developing an integer programming model and data pipeline with GIS mapping. GIS is a geospatial tool to present, display, and analyze data over space and time. We work with Asheville, North Carolina (NC) as a case study city and conduct extensive analyses to answer strategic decisions related to budget, equity, targeted demographics, and threshold metrics. The contributions of the paper are as follows:

- A location model to optimize the selection of new local parks to improve equitable access for residents.
- A person-centered modeling approach where access is parameterized according to key factors of resident experience.
- The presentation of geospatial data collection methods to support model inputs with multiple publicly available datasets.
- A case study of Asheville, NC that uses the park optimization model to consider policy questions and illustrate the potential for optimization in park planning.

The remainder of the paper is structured as follows. Section 2 introduces parks and recreation in Asheville, NC. Section 3 presents the formulation for the park equity model. Section 4 is a presentation of our data collection process. Section 5 provides the case study and model analyses. We discuss the results in Section 6 and conclude in Section 7.

**2. Parks and Recreation in Asheville, NC**

Asheville is a community located near the Appalachian Mountains in the western portion of NC. The city is home to approximately 94,000 people and has many local parks and a handful of national greenspaces (US Census Bureau, 2021b). The city manages 64 existing parks that are shown in Figure 1 (*Asheville Parks*, 2019). We will focus our analysis on the 52 that are freely available and open for general use (Section 4.3). The distribution of current parks is mainly focused within the central and eastern regions, whereas the northern, southern, and western regions have fewer. Only 3% of Asheville's 44.9 square-mile land area is composed of parks and greenspaces (The Trust for Public Land, 2022a; US Census Bureau, 2021b).

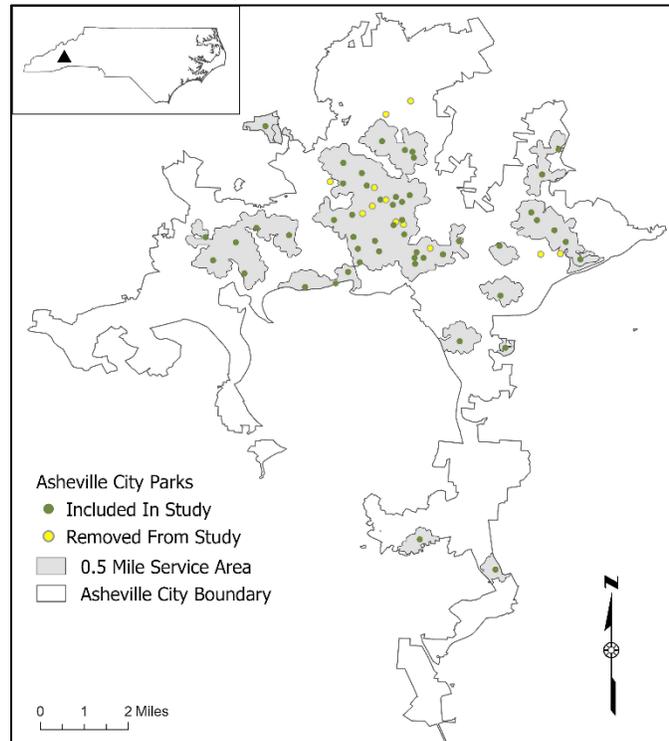

**Fig 1.** Distribution of parks across the city of Asheville

One measure of park access is distance from residents to parks. A standard threshold is to be within a 10-minute walk (The Trust for Public Land, 2022c). Across all Asheville residents, only 44% meet this criterion (The Trust for Public Land, 2022a). Access varies across racial-ethnic subpopulations. Just above half of Black residents live within a 10-minute walk (56%) whereas only a third of Hispanic residents do (34%), a range of 22% across racial-ethnic subpopulations. Forty-three percent of white residents are within a 10-minute walk.

We selected Asheville, NC as the subject of our case study because initiatives toward local park and greenspace equity already exist. The Asheville Parks and Recreation Department (APRD) created a Racial Equity Action Plan to combat injustices encountered by people in marginalized racial-ethnic groups (Asheville Parks and Recreation, 2022). Objectives include increasing environmental justice, support for community needs, and celebration of diversity. The department addresses these goals through representation, examination and prevention of the main sources of park inequities, and developments in policy. One aim is to allocate money and human resources to reduce inequities.

A current mechanism to provide funding is the Rehabilitation Community Investment Program. To identify target sites for improvements, the APRD uses a point-based method of neighborhood and park characteristics to represent the need for improvement in regions of the city (Asheville Parks and Recreation, 2022). It gives priority to areas with greater population density, poverty, crime, and numbers of youth and to parks that are worn, old, and unmaintained. The process is effective though time-consuming. Compiling and analyzing this information requires substantial resources that other cities may not possess. One opportunity to scale data collection and analysis is through the use of GIS-enabled optimization. Publicly available datasets for local, regional, and national organizations are readily available, and several cities, including Asheville, have a GIS department. To support efforts to improve park equity, we develop an optimization model that uses publicly available data on park and community composition. It can be used to support analyses for communities that may not have the resources to conduct large-scale studies of community need. We further use it to study policies in Asheville that may be difficult to evaluate with a point-based system.

**3. Park Equity Model**

In this section, we introduce the park equity model. We begin by presenting the definition and parameterization of park access (Section 3.1) and an overview of the model structure (Section 3.2). The formulation is in Section 3.3. Two variants of the model – an alternative objective function and uncapacitated approach – are presented in Section 3.4. Model assumptions are in Section 3.5, and structural properties are presented in Section 3.6.

**3.1 Park Access**

Environmental justice is multi-faceted, and representing access is difficult. In our context, we define *access* to parks to be the composition of four elements: (i) distance from residents to parks, (ii) appropriate capacity, (iii) heat, and (iv) tree cover. We define thresholds for each element; if a resident's primary park is within the thresholds, they experience good access. If a resident's park is outside of a threshold, it reflects a deviation from good access.

Consider an example of two resident locations that have different primary parks in Figure 2. The resident locations are labeled 1 and 2, and their primary parks are scored according to the four elements of access. The thresholds for good access are designated with brackets. Resident location 1 is further from its park than is ideal (outside the threshold) but has good access according to the other three metrics (within the thresholds). Resident location 2 is close to its park, and the park has acceptable tree cover (within the thresholds); however, it is over capacity and has excess heat (outside the thresholds).

These access elements were selected from a list of approximately two dozen factors. The initial list included factors related to distance, crowding, environment, land cover, transportation (e.g., bus stops), safety, amenities, and facilities. It was gathered from literature (Bruton & Floyd, 2014; Engelberg et al., 2016; Humpel et al., 2002; Kaczynski et al., 2020; National Recreation and Park Association, 2011; Suminski et al., 2012; Williams et al., 2020) and through discussions with parks and recreation planners and environmental justice experts. To focus the analysis on a parsimonious group of factors that could be parameterized with publicly available data and could influence strategic decisions in the case study city, we selected three categories of factors (distance, crowding, and environmental [comprised of tree cover and heat]). They were jointly decided across the co-authors representing operations researchers, parks and recreation planners, and environmental justice experts. The selection of model features was validated in meetings with the Director of the APRD and a manager of Urban Design Strategies. The elements of access are not intended to be exhaustive, and if the model is applied in other cities, different elements may be selected; generalizability of the model structure is further discussed in Section 6.3.

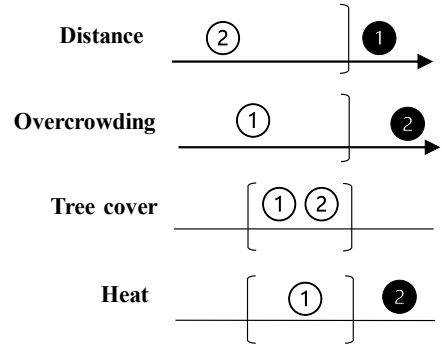

**Fig. 2.** Example values for elements of access for two resident locations with different primary parks

**3.2 Model Structure**

To address questions of strategic park planning, we propose a discrete location model of park location and access. The key sets are the resident locations, $L$; parks, $K$; and race/ethnicities, $R$. The core decisions are (i) where to locate parks and (ii) which park each resident location will primarily visit. A resident's access to parks is comprised of four elements that are characteristics of their primary park: walking distance, capacity, heat, and tree cover. If a resident's park does not fall within the acceptable ranges of an element, the deviation is recorded as a slack variable. The goal is to minimize weighted deviations from ideal access. We consider two styles of objective functions; the first minimizes weighted deviations for the subpopulation with the highest (Section 3.3), and the second minimizes weighted deviations across all residents (Section 3.4).

More specifically, the objective functions consider the extent to which each resident's primary park is out-of-bounds of ideal access across the four dimensions. This is a coverage-style approach. Deviations are multiplied by a weight, $w$, that indicates its relative priority vs. other categories and a normalization value, $n$, that converts each category to the same order of magnitude. These are further multiplied by the number of people $t_{lr}$ of each race/ethnicity $r \in R$ who live in each resident location $l \in L$. The granularity of the set of resident locations $L$ is based on the decision-maker and data availability; for example, it may represent neighborhoods or Census-defined areas. We include a strategic priority weight for each subpopulation $r \in R$, $q_r$, to be able to emphasize access for particular groups.

For example, consider resident location 1 in Figure 2. If the distance to the park is 0.75 miles and the threshold distance is 0.5 miles, then the distance deviation is 0.25 miles. This deviation value is multiplied by the weight on the distance element (indicating how important it is to be within the threshold), a normalization term (to convert the deviation to the same scale as the other elements), and the number of people who experience the deviation.

### 3.3 Notation and Formulation

Table 1 lists the sets, parameters, and decision variables, and the formulation is as follows.

**Table 1.** Sets, parameters, and decision variables of the park equity model formulation.

| | |
|---|---|
| ***Sets*** | |
| $K$ | set of all park locations ($K = K^{existing} \cup K^{candidate}$) |
| $K^{existing}$ | set of existing parks |
| $K^{candidate}$ | set of candidate parks |
| $L$ | set of resident locations |
| $R$ | set of race/ethnicities |
| ***Parameters*** | |
| $e_k$ | $:= \begin{cases} 1 & \text{if park already exists at park location } k \in K \\ 0 & \text{otherwise} \end{cases}$ |
| $d_{kl}$ | walking distance from resident location $l \in L$ to park location $k \in K$ |
| $m$ | maximum allowable distance for any resident to their primary park |
| $a_k$ | capacity of park location $k \in K$ |
| $b$ | budget to purchase park locations |
| $f_k$ | cost to purchase park location $k \in K$ |
| $c_k^+, c_k^-$ | amount of heat above, below desirable range for park location $k \in K$ |
| $v_k^+, v_k^-$ | amount of tree cover above, below desireable range for park location $k \in K$ |
| $t_{lr}$ | number of people of demographic group $r \in R$ who live in location $l \in L$ |
| $w^{dist+}$ | penalty weight of excess distance |
| $w^{cap+}$ | penalty weight of park over capacity |
| $w^{heat+}, w^{heat-}$ | penalty weight of excess, deficit park heat |
| $w^{tree+}, w^{tree-}$ | penalty weight of excess, deficit park tree cover |
| $q_r$ | strategic emphasis of demographic group $r \in R$ |
| $n^{dist}$ | normalization of distance |
| $n^{cap}$ | normalization of capacity |
| $n^{heat}$ | normalization of heat |
| $n^{tree}$ | normalization of tree cover |
| ***Decision Variables*** | |
| $y_k$ | $:= \begin{cases} 1 & \text{if park location } k \in K \text{ is selected} \\ 0 & \text{otherwise} \end{cases}$ |
| $x_{kl}$ | $:= \begin{cases} 1 & \text{if } k \in K \text{ is the primary park location of residents in resident location } l \in L \\ 0 & \text{otherwise} \end{cases}$ |
| $\alpha_r$ | total weighted deviations for demographic group $r \in R$ |
| $\alpha^{max}$ | total weighted deviations for demographic group with the highest |
| $\tilde{d}_l^+$ | extra distance to primary park for resident location $l \in L$ |
| $\tilde{a}_k^+$ | amount that park $k \in K$ is over capacity |
| $u_l^{dist}$ | $:= \begin{cases} 1 & \text{if distance to primary park is within } m \text{ for location } l \in L \\ 0 & \text{otherwise} \end{cases}$ |
| $u_k^{cap}$ | $:= \begin{cases} 1 & \text{if park } k \in K \text{ is within capacity, } a_k \\ 0 & \text{otherwise} \end{cases}$ |

$$\text{minimize } \alpha^{max} \tag{1}$$

*Subject to:*

$$\alpha_r = \sum_{l \in L} \sum_{k \in K} \left( q_r t_{lr} \left[ n^{dist} w^{dist+} \tilde{d}_l^+ + \begin{pmatrix} n^{cap} w^{cap+} \tilde{a}_k^+ \\ + n^{heat} w^{heat+} c_k^+ \\ + n^{heat} w^{heat-} c_k^- \\ + n^{tree} w^{tree+} v_k^+ \\ + n^{tree} w^{tree-} v_k^- \end{pmatrix} x_{kl} \right] \right) \quad \forall r \in R \tag{2}$$

$$\alpha^{max} \geq \alpha_r \quad \forall r \in R \tag{3}$$

$$\sum_{k \in K} x_{kl} = 1 \quad \forall l \in L \tag{4}$$

$$x_{kl} \leq y_k \quad \forall k \in K, l \in L \tag{5}$$

$$e_k \leq y_k \quad \forall k \in K \tag{6}$$

$$\sum_{k \in K} f_k y_k \leq b \tag{7}$$

$$\sum_{k \in K} d_{kl} x_{kl} - \tilde{d}_l^+ \leq m \quad \forall l \in L \tag{8}$$

$$\tilde{d}_l^+ - (1 - u_l^{dist}) \left( \sum_{k \in K} d_{kl} x_{kl} - m \right) \leq 0 \quad \forall l \in L \tag{9}$$

$$\sum_{l \in L} \sum_{r \in R} t_{lr} x_{kl} - \tilde{a}_k^+ \leq a_k \quad \forall k \in K \tag{10}$$

$$\tilde{a}_k^+ - (1 - u_k^{cap}) \left( \sum_{l \in L} \sum_{r \in R} t_{lr} x_{kl} - a_k \right) \leq 0 \quad \forall k \in K \tag{11}$$

$$y_k \in \{0,1\} \quad \forall k \in K \tag{12}$$

$$x_{kl} \in \{0,1\} \quad \forall k \in K, l \in L \tag{13}$$

$$\alpha^{max} \geq 0 \tag{14}$$

$$\alpha_r \geq 0 \quad \forall r \in R \tag{15}$$

$$\tilde{d}_l^+ \geq 0 \quad \forall l \in L \tag{16}$$

$$\tilde{a}_k^+ \geq 0 \quad \forall k \in K \tag{17}$$

$$u_l^{dist} \in \{0,1\} \quad \forall l \in L \tag{18}$$

$$u_k^{cap} \in \{0,1\} \quad \forall k \in K \tag{19}$$

The objective function (1) minimizes the maximum weighted deviations, i.e., the weighted deviations of the subpopulation with the highest. Constraints (2) calculate the total weighted demographic deviation for each subpopulation. Total deviations are a function of the normalized and weighted values of walking distance, capacity, heat, and tree cover. An exact linearization of this constraint is provided in Appendix A (22-26). Constraints (3) determine the maximum weighted demographic deviations across all subpopulations. Constraints (4) ensure that all resident locations have a primary park. Constraints (5) state that residents may only visit selected parks. Constraints (6) require that a park be selected if it already exists. Constraint (7) limits spending to the available budget.

Constraints (8) require residents to be within the desirable distance of their primary park, with slack variables to record deviations. Constraints (9), in combination with constraints (8), enforce an exact value for the slack variable $\tilde{d}_l^+$ for each resident location $l \in L$. Constraints (9) are nonlinear, and the exact, implemented, linearizations are provided in Appendix A (constraints (27-31)). Note that without constraints (9), the slack variables may have values that are artificially higher than the actual distance deviation to improve "equity" between groups. Similarly, constraints (10) require that the number of residents who primarily visit a park be within the park's capacity, with slack to record deviations. Constraints (11) ensure that these represent the true values over capacity. Exact linearizations, that were implemented in place of nonlinear (11), are provided in Appendix A (constraints (32-36)). Constraints (12-19) are domain constraints.

### 3.4 Model Variants

We consider two variants of the base model. One is an alternative objective function (20) that minimizes overall weighted deviations, across all subpopulations.

$$\text{minimize} \sum_{r \in R} \alpha_r \tag{20}$$

A second variant removes capacity limits for each park location. In this model, either objective function (1) or (20) can be used. The constraints include (3)-(9), (12)-(13), (16), and (18), which are defined above, and one additional set of constraints (21).

$$\alpha_r = \sum_{l \in L} \sum_{k \in K} \left( q_r t_{lr} \left[ n^{dist} w^{dist+} \tilde{d}_l^+ + \begin{pmatrix} n^{heat} w^{heat+} c_k^+ \\ + n^{heat} w^{heat-} c_k^- \\ + n^{tree} w^{tree+} v_k^+ \\ + n^{tree} w^{tree-} v_k^- \end{pmatrix} x_{kl} \right] \right) \quad \forall r \in R \tag{21}$$

### 3.5 Model Assumptions

The goal of this paper is to develop a modeling framework, compile datasets, and illustrate the potential for optimization to address questions of park access and recreation planning. To support these goals, the models include the following simplifying assumptions, developed in partnership with park planning experts and stakeholders.

The first is that the model considers a single, aggregated measure of access, $\alpha_r$ for each demographic group $r \in R$. An alternative would be to consider disaggregated elements of access and use a multi-objective approach. In this paper, we use an aggregate measure to make it easier for decision-makers to interpret the results. One downside to aggregation is that there is the risk that at "optimality" residents may experience disparities in lower-weighted elements of access. To evaluate if this occurs, we present disaggregated, unweighted results (Figure 8). If these results are undesirable, the model could be re-run with different priority weights.

We only consider a resident's primary park to evaluate park access. This is consistent with the park planning community's emphasis on providing access to at least one park for all residents (The Trust for Public Land, 2022a), though in practice, a resident may visit any number of parks. Actual access may be higher than the model would suggest, and this underestimation of access is particularly true when the set of resident locations $L$ represents large areas, as in our case study (Section 4.1). We note we parameterize $L$ by Census block groups because this is the most granular location that has race/ethnicity data. One substantial assumption is that the capacitated models define primary parks through the lens of the municipal planner with best-case access, rather than a resident-choice approach. The model is structured in this way due to lack of data on park visitation but was considered a reasonable approximation by stakeholders. To mitigate this limitation, throughout, we compare the capacitated models with uncapacitated models. The latter do not have this assumption because coordination between the residents does not occur without capacity constraints. We would expect resident decisions in practice to be between the capacitated and uncapacitated versions.

The models are deterministic, including deterministic population sizes. This is consistent with seminal descriptive analyses on park access (Ibes, 2015; Talen & Anselin, 1998) and state-of-the-art non-optimized recommendations for new locations (The Trust for Public Land, 2022a). This modeling choice facilitates ease of comparison with existing tools as well as ease of understanding by stakeholders due to conceptual familiarity. We assume that residents in each location $l \in L$ have the same primary park. This is consistent with the framing of a *neighborhood* park where residents in the same area view a particular park as theirs, per discussions with park planners. This could be relaxed by further stratifying $x_{kl}$. We do not directly consider park sizing decisions, though differing parcel sizes in the input data implicitly address this feature in part. To consider a limited number of multi-parcel parks, the modeler could define new mutually-exclusive park options $k \in K$ that are comprised of multiple parcels within the input data, cf. facility location sizing (Correia and Captivo, 2003).

### 3.6 Structural Properties

An analysis of the base model leads to the following observations. Proposition 1 states that insufficient park access may be inherent to certain city structures regardless of optimization or investment. For example, cities with limited land to use for parks may be unable to satisfy targets for acres per person. We also observe that as parks are added (Proposition 2), $\alpha^{max}$ remains the same or improves. This suggests that as extra land becomes available in underserved areas, it may be worthwhile to select it as a park. Similarly, budget increases (Corollary) may lead to the same or improved worst-case access. If extra money is allocated or donated for park purchases, this would not worsen access, as defined by the model.

We define additional notation as follows. Let $y_k^1, y_k^2$ indicate two selections of parks $k \in K$ feasible to (2)-(19), where $y_k^1$ is nested within $y_k^2$. That is, $y_k^1 = y_k^2, \forall k \in K \setminus K'; y_{k'}^1 = 0, y_{k'}^2 = 1 \forall k' \in K'$. Let $\alpha^{max}(\cdot)$ represent the objective value (1) for a problem instance with a given solution or parameter as its argument.

**Proposition 1.** If $\sum_{k \in K} a_k < \sum_{l \in L} \sum_{r \in R} t_{lr}$ or $\exists l \in L: \min_k d_{kl} > m$, then $\alpha^{max} > 0$.

**Proof.** $\sum_{k \in K} a_k < \sum_{l \in L} \sum_{r \in R} t_{lr}$ implies $\sum_{k \in K} \tilde{\alpha}_k^+ > 0$ by Constraints (4), (10), and (13). $\exists l \in L: \min_k d_{kl} > m$ implies $\sum_{l \in L} \tilde{d}_l^+ > 0$ by Constraints (4) and (8). If either $\sum_{k \in K} \tilde{\alpha}_k^+ > 0$ or $\sum_{l \in L} \tilde{d}_l^+ > 0$, then $\alpha^{max} > 0$ by Constraints (2)-(4).

**Proposition 2.** Given nested solutions $y_k^1, y_k^2$, then $\alpha^{max}(y_k^1) \geq \alpha^{max}(y_k^2)$.

**Proof.** Access $\alpha_r$ for group $r \in R$ is a function of the resident locations $l \in L$ associated with each park $k \in K$, $x_{kl}$, per (2), (4), (9), and (10). All $x_{kl} \forall k \in K, l \in L$ that are feasible with $y_k^1$ are also feasible with $y_k^2$, per (5). For parks $k' \in K'$, any resident location $l \in L$ such that $x_{k'l} = 1$ that are feasible with $y_{k'}^2$ are infeasible with $y_{k'}^1$, per (5). This implies $\alpha^{max}(y_k^1) \geq \alpha^{max}(y_k^2)$ per (3) and (1).

**Corollary.** Given $b^2 > b^1$, $\alpha^{max}(b^1) \geq \alpha^{max}(b^2)$.

**Proof.** Let $y_k^1, \forall k \in K$ be feasible to (2)-(19) with budget $b^1$. By (7), $y_k^1$ is feasible with budget $b^2$, implying $\alpha^{max}(b^1) \not< \alpha^{max}(b^2)$.

## 4. Data Collection and Model Inputs

We conducted extensive data collection to translate data from publicly available databases to model inputs. In this section, we present the approach for our case study of Asheville, NC. The set of resident locations is defined in Section 4.1, and the population sizes for each demographic group in each location are derived in Section 4.2. We define the set of existing park locations and candidate park locations in Section 4.3 and 4.4, respectively. Land costs are estimated in Section 4.5. The capacities of each existing and candidate park are provided in Section 4.6 and environmental factors in Section 4.7. Walking distances of each resident location to each candidate park are presented in Section 4.8. Finally, parameters to support the case study analyses are presented in Section 4.9. A compilation of each parameter with its source information is presented in Appendix Table B1. It includes descriptions, GIS data type, year, and source for each set or parameter.

### 4.1 Resident Locations and Race/Ethnicities ($L, R$)

We define the set of resident locations, $L$, as the 2019 Census block groups that are within Asheville city limits (*Asheville City Limits*, 2017; US Census Bureau, 2021a). This is the smallest geographic unit that has race/ethnicity data available at the time the analysis was conducted. To focus on residents within the city's jurisdiction, where block groups may extend beyond Asheville's city limits, we clip the geographic areas and only include the portion of block groups that are within city limits (Appendix Figure B1). We also limit the block groups to those which have at least 25 people (*Race (P1)*, 2020), eliminating 11 of the original 88 block groups.

To define the set of race/ethnicities, we followed the US Census Bureau's six classifications of White, Black or African American, American Indian and Alaska Native, Asian, Native Hawaiian and Other Pacific Islander, and Some Other Race (*Race (P1)*, 2020), with slightly abbreviated labels, $R =$ {White, Black, Indigenous, Asian, Pacific Islander, Other}.

### 4.2 Population Count by Race/Ethnicity and Location ($t_{lr}$)

The number of people of each race/ethnicity $r \in R$ that live in resident location $l \in L$ is defined as $t_{lr}$. The Census race data contains the number of people who live in each Block Group in 2020 (*Race (P1)*, 2020). In this data, each individual may either belong to a single or multiple race/ethnicities. The population size of each race/ethnicity $r \in R$ is the number of people of a single race/ethnicity in addition to people who list the race/ethnicity under multiple race associations. Note that this grouping method double-counts individuals with multiple race/ethnicities; this approach was selected because the primary objective (1) considers subpopulations individually, and it includes people with multiple races/ethnicities in each relevant group.

The population count data is available for 2020 Block Groups (*Race (P1)*, 2020), but the resident locations $l \in L$ in the model are the 2019 Block Groups (US Census Bureau, 2021a). Because of redistricting, a few of the areas change from 2019 to 2020 (Appendix Figure B2). To determine population count $t_{lr}$ by race/ethnicity $r \in R$ for each location $l \in L$, we applied the 2020 Block Group data to the 2019 Block Groups proportionally to the percentage of overlap.

If location $l \in L$ had a portion of its area removed because it was outside of city boundaries (Appendix Figure B1), we proportionally reduced the population sizes by the amount of area that was removed. For example, if 50% of a block group's area was within city limits, the population size used in the model is 50% of the population in the block group. This assumes homogeneity of population density. The number of people in each resident location $l \in L$ is presented in Figure 3 ($\sum_{r \in R} t_{lr}$).

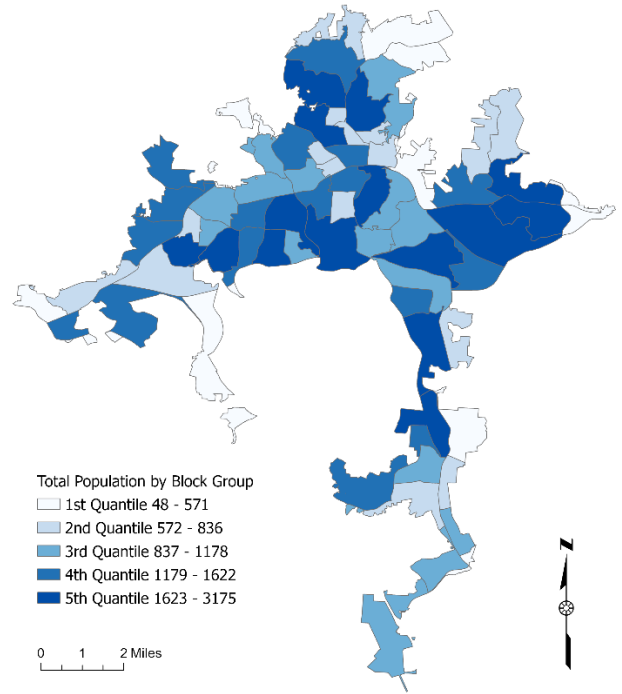

**Fig. 3.** Population size in each resident location $l \in L$

### 4.3 Locations of Existing Parks ($K^{existing}$)

For the purpose of these analyses, we will only include parks that are maintained by the city, open for general use, and freely available, consistent with the definition from the Trust for Public Land (The Trust for Public Land, 2022a). The City of Asheville lists 64 parks (Figure 1) (*Asheville Parks*, 2019). We checked each park to determine whether it satisfies the inclusion criteria and removed twelve parks from the initial list (City of Asheville, 2022; Google, 2021). Reasons for park removal are payment required for entry, single purpose or professional team use only, and limited infrastructure or space. We do not include parks outside of city limits, and we expect edge effects to be minimal due to the 0.5 mile distance threshold for ideal access.

The final set of 52 existing parks is illustrated in Figure 1. The service area of each is defined as within a 0.5-mile walking distance. We note that there are several areas within Asheville that are not currently covered by an existing city park. The map also depicts the listed parks that did not meet the inclusion criteria in yellow. Some of these are within the service areas of existing parks, indicating that residents are still covered by a park.

### 4.4 Candidate Parks ($K^{candidate}$)

To define a set of candidate parks, we considered single parcels within Asheville that meet certain characteristics. Parcels were collected via the shapefile of Buncombe County parcels (Buncombe County GIS, 2023). Inclusion criteria included land compactness, size of at least one acre, and not located within a flood zone. We also only selected candidate sites within building zones that allow greenspace development; this restriction eliminated the airport zone (*Municode: Asheville, NC*, 2022). Sites should not be within protected land areas nor service areas of existing parks. For the latter, we removed parcels located within a 0.5-mile walking distance from the existing parks ($K^{existing}$). We verified that selected parcels were not roadways, parking lots, or land for existing structures via ArcGIS basemap visual satellite imagery; this process ensures that no housing demolition occurs to create park space. There were a total of 39,480 parcels within Asheville, and the selection process led to a set of 138 candidate new park locations, $K^{candidate}$. The set of candidate parks is shown in Figure 4 along with elements of existing parks and service areas, water features, floodways, and protected land.

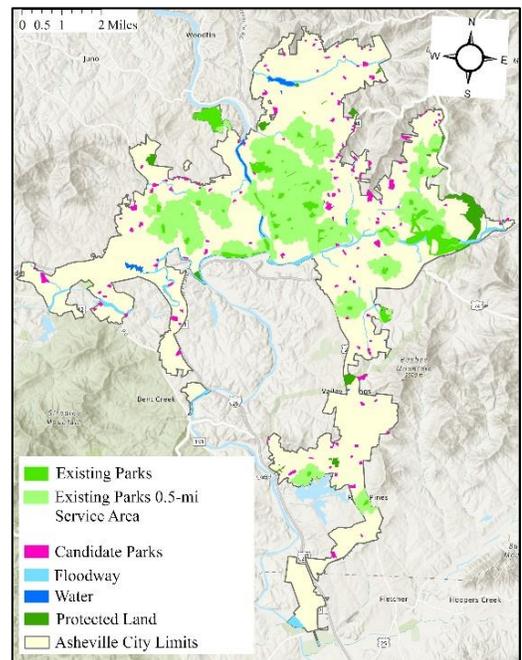

**Fig. 4.** Candidate parks, existing parks, and geographic elements.

### 4.5 Park Cost ($f_k$)

We estimated the cost to purchase each park location $k \in K$ via available land values. Land values were available in the Buncombe County parcels attribute table (Buncombe County GIS, 2023). Any parcels already owned by the City of Asheville, e.g., $k \in K^{existing}$, had no cost.

Some of the parcels did not have land values listed. To estimate the cost in these cases, we multiplied the number of acres in the candidate park (ESRI, 2023) by the approximate cost-per-acre. To determine the approximate cost-per-acre, we divided the city into zones and calculated the average cost-per-acre for parcels with a land value listed. The distribution of candidate parks within each cost zone as well as whether each candidate park had a listed land value is provided in Appendix Figure B3.

### 4.6 Park Capacity ($a_k$)

Park capacity, $a_k$, is the number of residents that park location $k \in K$ can accommodate. Capacity varies by park size, and there should be at least one acre of park land for every 100 residents that visit a park (Moeller, 1965). Candidate park sizes are set as the size of the parcel. To calculate capacity, we multiplied the acres of each park location $k \in K$ (ESRI, 2023) by 100 and rounded down.

### 4.7 Park Environmental Factors ($c_k^+, c_k^-, v_k^+, v_k^-$)

Two key environmental factors are park heat and tree cover. Both have a range of allowable values, and deviations are penalized in the model. For each park, we calculated the average heat and tree cover and compared them to the allowable ranges. For candidate parks, the values represent the heat and tree cover of the parcels eligible to be purchased. The resulting excess and deficit park heat ($c_k^+$ and $c_k^-$) and excess and deficit park tree cover ($v_k^+$ and $v_k^-$) are then included in the model as measures of deviation. Note that these parameters effectively represent slack variables if park location $k \in K$ is selected; they are pre-calculated to avoid non-linearity.

The data sources for heat index (The Trust for Public Land, 2022b) and tree cover (Multi-Resolution Land Characteristics Consortium, 2020) provided raster datasets. With raster data, the data is presented in "cells." Each cell represents a region (approximately 30m x 30m in the included datasets), and the value of the cell is the majority value within this region.

The heat index dataset presented values for heat relative to Asheville's average heat on a scale of 1 to 5. An index of 1 represented a cell with heat slightly above the city average while an index of 5 represented a heat greatly larger than the city average. Cells without a heat index value experience heat less than or equal to the city average. Tree cover was presented as a percentage of vegetation within each raster cell (0%-100%). The park locations, $K$, may include multiple cells and cells that overlap with the park boundaries. To determine heat and tree cover for the park location $k \in K$ overall, we resampled the cells to produce smaller cell sizes of 5m x 5m. Any smaller cell that was completely within the boundaries of park location $k \in K$ was used to calculate an average heat index and tree cover for park $k$ (ESRI, 2023). The averages were compared to the allowable range, and any deviations were recorded as $(c_k^+, c_k^-)$ and $(v_k^+, v_k^-)$ for heat and tree cover of park location $k \in K$, respectively.

## 4.8 Walking Distance ($d_{kl}$)

The walking distance from resident location $l \in L$ to park $k \in K$ is given by $d_{kl}$. To calculate the distances, we used the Network Analysis Toolbox in ArcGIS Pro (ESRI, 2023). The origins (resident locations $l \in L$) and destinations (parks $k \in K$) were represented by their geometric center-points. The distance was calculated via the shortest path between each origin and destination. The eligible routes were roadways that allow pedestrian travel within the ArcGIS Pro software (ESRI, 2023).

## 4.9 Parameters for Analyses

Additional parameters include coefficients for the objective functions and thresholds for the allowable ranges. Each of these can be adjusted based on user priorities. Table 2 provides the selected values of these parameters. For consistency, we hold these parameters constant throughout all model analyses. We place the greatest amount of importance upon distance as a category of access. We aim for all residents to be within 0.5 miles of a park and to have at most 100 residents per park acre (Moeller, 1965). We determined the acceptable heat range to target areas that experience relatively high amounts of heat. We determined the acceptable tree cover range such that we target areas that have a moderate amount of vegetation. This enables the model to select park sites that may support the development of multiple amenities while maintaining the provision of shade. We normalized access categories by projecting them onto the same numerical scale relative to the maximum value for each parameter. The baseline strategic emphasis parameter for each subpopulation is set to 1 for all analyses, except the prioritization analysis in Section 5.3.

**Table 2.** Objective function coefficient values

| Parameter | Weight ($w$) | Normalization ($n$) | Threshold |
|---|---|---|---|
| Distance | 0.9 | 5 | |
| Capacity | 0.25 | 1/150 | |
| Excess heat | 0.2 | 20 | 4 |
| Deficit heat | 0.05 | 20 | 1 |
| Excess tree cover | 0.25 | 1 | 70 |
| Deficit tree cover | 0.2 | 1 | 20 |

## 5. Case Study Results

Using the city of Asheville, NC as a case study, we conduct analyses to address policy questions of budget use, strategic targeting, and metrics of access. Specifically, these are: (1) how does increasing the budget for park investments affect access to parks? (2) What is the impact of long-term vs. myopic park planning? (3) How does strategically targeting demographic groups affect park plans? (4) How does changing an access threshold affect park selection? These are addressed in Sections 5.1-5.4, respectively. In Section 5.5, we compare access-oriented outcomes to those from classic location-allocation models.

In each of the policy analyses, we consider multiple model instances: two types of objective functions (minimizing the maximum weighted deviations [objective function (OF) 1] and minimizing overall weighted deviations [OF 20]) and two types of capacity (capacitated and uncapacitated). We designate these four combinations *Min Max Dev Cap* (OF 1, capacitated), *Min All Dev Cap* (OF 20, capacitated), *Min Max Dev Uncap* (OF 1, uncapacitated), and *Min All Dev Uncap* (OF 20, uncapacitated). The model is programmed in AMPL with Gurobi as the optimization solver (Fourer, Gay, and Kernighan, 2022; *Manual, Gurobi Opimizer Reference*, 2022). The analyses were conducted on a laptop with an Intel Core i7 processor, 2.5 GHz CPU, and 8 GB of memory. Each model instance is solved to exact optimality. Solution times are at most 5 minutes with OF 1 and at most 10 minutes with OF 20. The code and data are available on GitHub (Young et al., 2022).

**5.1 Budget Available**

Governmental and recreational organizations may allocate money to purchase new parks. In this subsection, we study the relationship between the budget to purchase new parks, $b$, and access. For each analysis, we consider a budget range of $0 to $3 million and analyze results at each increment of $250,000. The budget of $0 reflects the parks currently open in Asheville, NC and parcels currently owned by the city (Appendix Figure B4). The upper bound of $3 million was selected to be beyond the threshold at which the incremental benefit of additional budget is minimal. We evaluate weighted deviations across all demographic groups (total and by category), weighted deviations for the subpopulation with the worst access (total and by category), and unweighted deviations (by category).

First, we observe how budget affects overall deviations (Figure 5). As the budget increases, weighted deviations decrease; this means that it is possible to improve park access by strategically selecting locations for new parks. At $1 million, weighted deviations are 37-51% lower than $b = $0, and at $2 million, weighted deviations are 44-55% lower.

There are diminishing returns for budget increases. An increase from $0 to $250,000 decreases overall weighted deviations by 21% vs. a decrease of 10% for a budget increase from $250,000 to $500,000. That is, as access improves, it becomes increasingly difficult to address through strategic park selection. There is a point at which additional funding no longer improves access (approximately $4.75 million for *Min All Dev Cap* and $3.5 million for *Min All Dev Uncap*).

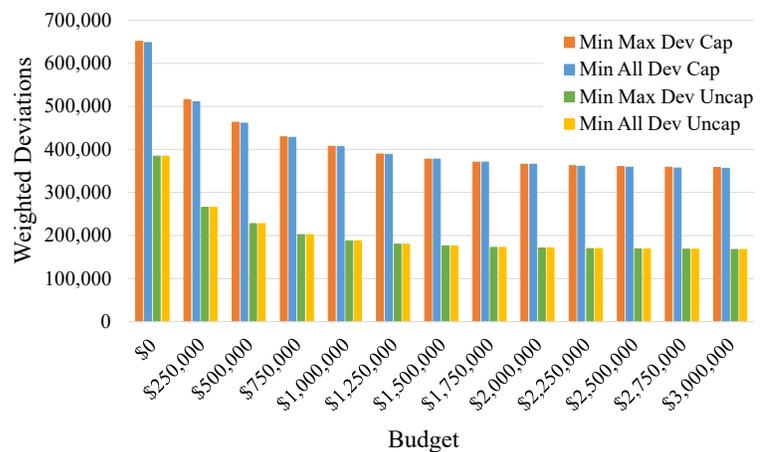

**Fig. 5.** Weighted deviations by budget

We also observe the growing size of park service areas with different budget amounts (Figure 6). With a budget of $0.5 million, parks are added around the periphery, with new locations in the southern, northern, and western regions of Asheville. With a budget of $1.5 million, parks are also added in the central region and in moderate- to high-population areas (Figure 3) such that the new service areas are contiguous with service areas of existing parks. The latter has the double benefit of reducing distance and capacity deviations at the existing parks.

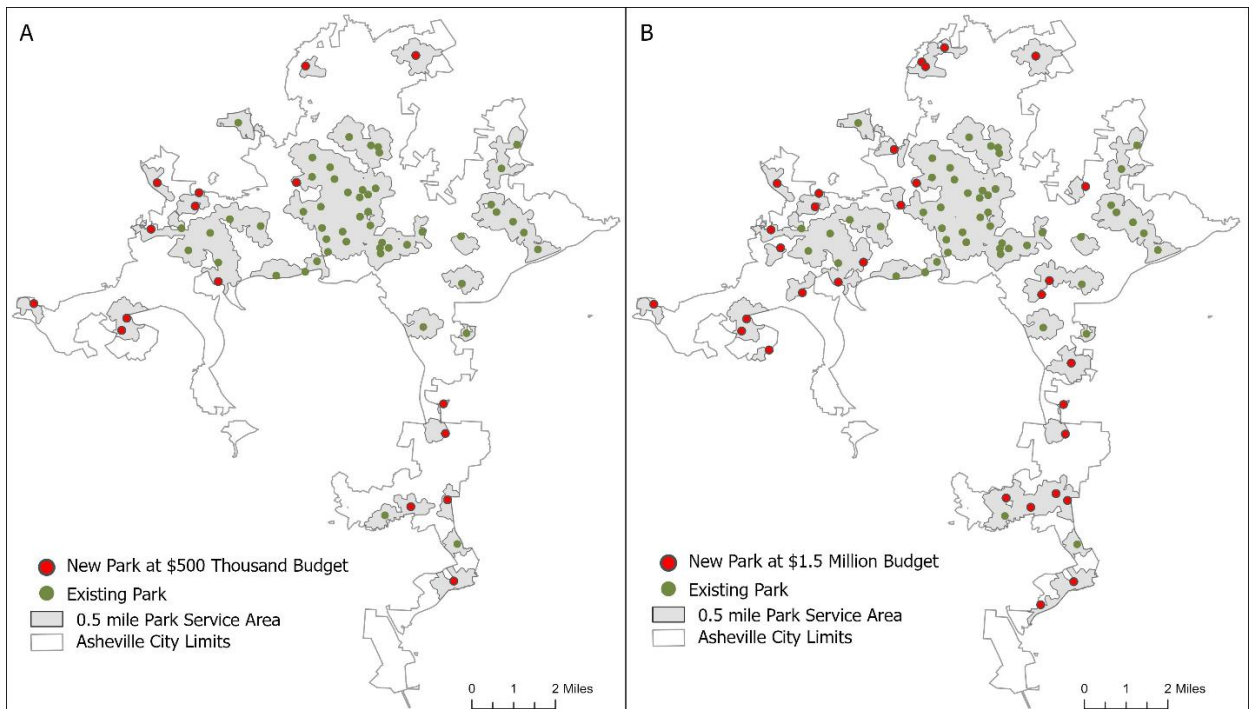

**Fig 6.** Selected park locations and area within 0.5 miles for budgets of $500,000 and $1,500,000 (*Min Max Dev Cap*)

Both objective functions (minimizing the deviations for the subpopulation with the highest and minimizing overall deviations) produce similar trends as the budget increases, and the resulting differences in overall weighted deviations are negligible for each budget amount. This indicates that overall access improves with both types of objective functions (*Min All* and *Min Max*). The two model types (capacitated and uncapacitated) lead to substantially different values for weighted deviations. Across all budget amounts, weighted deviations are 69-117% higher in the capacitated model. The capacitated models have higher weighted deviations because they penalize parks' overcapacity, whereas the uncapacitated models do not. As indicated by Proposition 1, we observe that there is a limit to which insufficient access can be improved given park locations $K$ in Asheville. We note that the lower bound on weighted deviations is reached at a smaller budget amount for the uncapacitated models than for capacitated. The uncapacitated models are not bounded by park capacity limits.

When we focus on the demographic subpopulation with the highest weighted deviations, $\alpha^{max}$, the results are similar (Appendix Figure C1). As budget increases, the weighted deviations decrease. With a budget of $1 million, the maximum weighted deviations are 38-52% lower than a baseline of $b = \$0$ (across capacitated and uncapacitated models). This suggests that strategically selecting parks can make a substantial impact in addressing limited park access. There are diminishing returns; at the first increment of $250,000, maximum weighted deviations decrease by 22-32% ($250,000 vs. $0), and at the second increment, maximum weighted deviations decrease by 10-15% ($500,000 vs. $250,000).

Up to this point, we have looked at aggregated measures of access, i.e., total weighted deviations across categories. We will now shift to consider the *composition* of the total deviations by evaluating budget effects on the disaggregated categories of access (distance, capacity, heat, and tree cover). We display results from the *Min Max Dev Cap* model (Figure 7). Note that the sum of the weighted deviations by category equals the total weighted deviations (e.g., Figure 5).

At the baseline budget of $0, there are weighted deviations in each of the four categories. Distance has the highest weighted deviations (42% of total), followed by capacity (31%), tree cover (17%), and heat (10%). As budget increases, weighted distance deviations decrease the most. It also becomes a smaller proportion; from 42% at baseline, it drops to 36% at $250,000 and 33% at a budget of $500,000. Capacity becomes a larger share of proportional weighted deviations (with 31% at baseline to 35% and 37% at $250,000 and $500,000, respectively). At some budget levels, decreasing distance comes at the expense of other categories. For example, weighted capacity deviations do not decrease monotonically; from $2 million to $2.25 million, weighted capacity deviations increase from 152,513 to 154,927 (relative increase of 1%). This results because the objective function optimizes weighted deviations across all categories rather than each individual factor. The proportion of weighted deviations by category continues to adjust as budget increases; at $1 million, the proportional contribution by category is (29%, 41%, 14%, 16%), and at $2 million it is (24%, 42%, 16%, 19%) for (distance, capacity, tree cover, and heat). Results for the maximum subpopulation weighted deviations are similar.

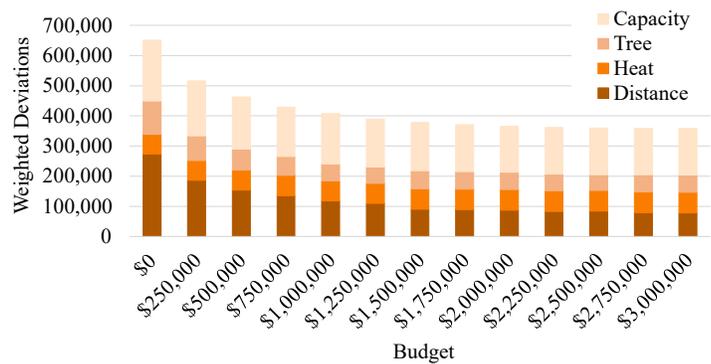

**Fig. 7.** Weighted deviations by category with different budgets (*Min Max Dev Cap* model)

Next, we consider *unweighted* distance deviations that represent the actual value of the slack variable, i.e., excess distance (Figure 8). These values are without the weighting and normalization terms used in the objective function. As the budget increases, the average and maximum distance deviations tend to decrease. With a budget of $1 million, parks are located such that residents are 0.4 miles further than ideal on average. Doubling the budget to $2 million leads to a set of parks in which residents are 0.3 miles further on average. The maximum resident-experienced distance deviation values are substantially higher than the average distance deviations. This suggests that several resident locations are nearly at a desirable distance from their primary park, but there are extremes at which this is not the case.

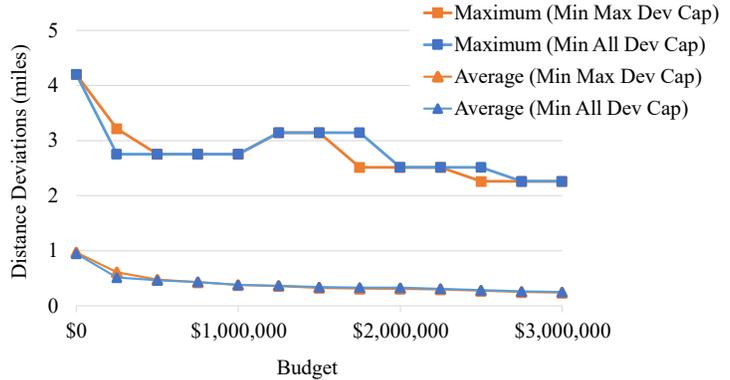

**Fig. 8.** Resident-experienced distance deviations

In the cases where distance deviations increase with budget, another category of access decreases (e.g., capacity in Appendix Figure C2). There is a pronounced trade-off between reducing distance and capacity deviations, and in the results from the uncapacitated models, there is a monotonic decrease in maximum distance deviation as budget increases (i.e., when capacity deviations are excluded from the objective). As budget increases, we note that the average resident-experienced capacity deviations tend to decrease as well. The average capacity overage is 1,266 people at baseline; 941 with $1 million budget; and 773 with $2 million budget (*Min Max Dev Cap*).

**5.2 Long-term vs. Myopic Planning**

Park selection decisions are affected not only by budget amount but also by when plans are made. To evaluate the potential effects of optimization during strategic planning, we consider two scenarios over a ten-year time horizon. In the long-term scenario, the full purchasing budget is known at the beginning of park planning, and the budget in the model, $b$, is set equal to the full 10-year budget. The actual purchase may be made over time as the money becomes available, but the locations would be pre-decided. In the myopic scenario, an allocation is made each year, and parks are purchased ad hoc. The budget, $b$, is set equal to the one-tenth of the 10-year budget; selected parks are added to the set of existing parks, $K^{existing}$ for any subsequent years. Unspent funds are not carried over. Note that in both scenarios, the same total amount of money may be allocated. We evaluate both scenarios at two different ten-year budget levels ($1 million and $2.5 million) using the model type *Min Max Dev Cap*. Outcomes are weighted deviations and locations.

For both total budget amounts, access is better with long-term planning (Figure 9). Overall weighted deviations are 5% lower with long-term planning with a budget of $1 million and 1% lower with a budget of $2.5 million vs. myopic planning. The benefit of long-term planning decreases as the budget increases. Similar differences in weighted deviations occur for the demographic group with the highest; weighted deviations are 5% and 1% lower vs. myopic planning for this group with $1 million and $2.5 million budgets, respectively. Both long-term and myopic planning improve access vs. status quo. Long-term planning could lead to better outcomes than myopic, but either would be an improvement; this validates the convention to use money to purchase parks when it is available.

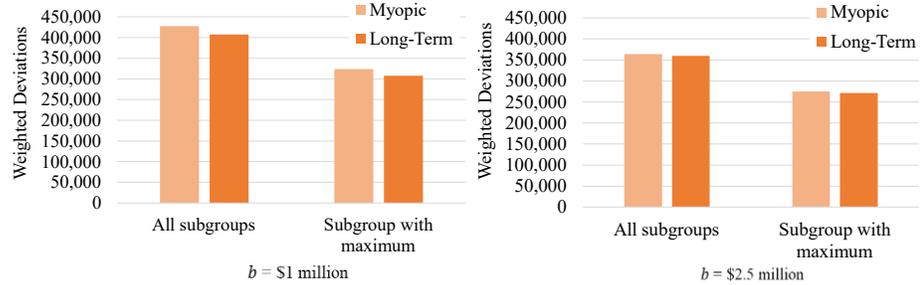

**Fig. 9.** Weighted deviations resulting from long-term vs. myopic planning

Next, we look at which parks are selected in each planning scenario. Figure 10 presents both scenarios at a budget level of $1 million (panel A) and at $2.5 million (panel B). The figures show Asheville's existing parks, the myopic (iterative) purchases each year, and the parks selected with the long-term planning budget. With a $1 million budget (Figure 10A), the model primarily selects candidate

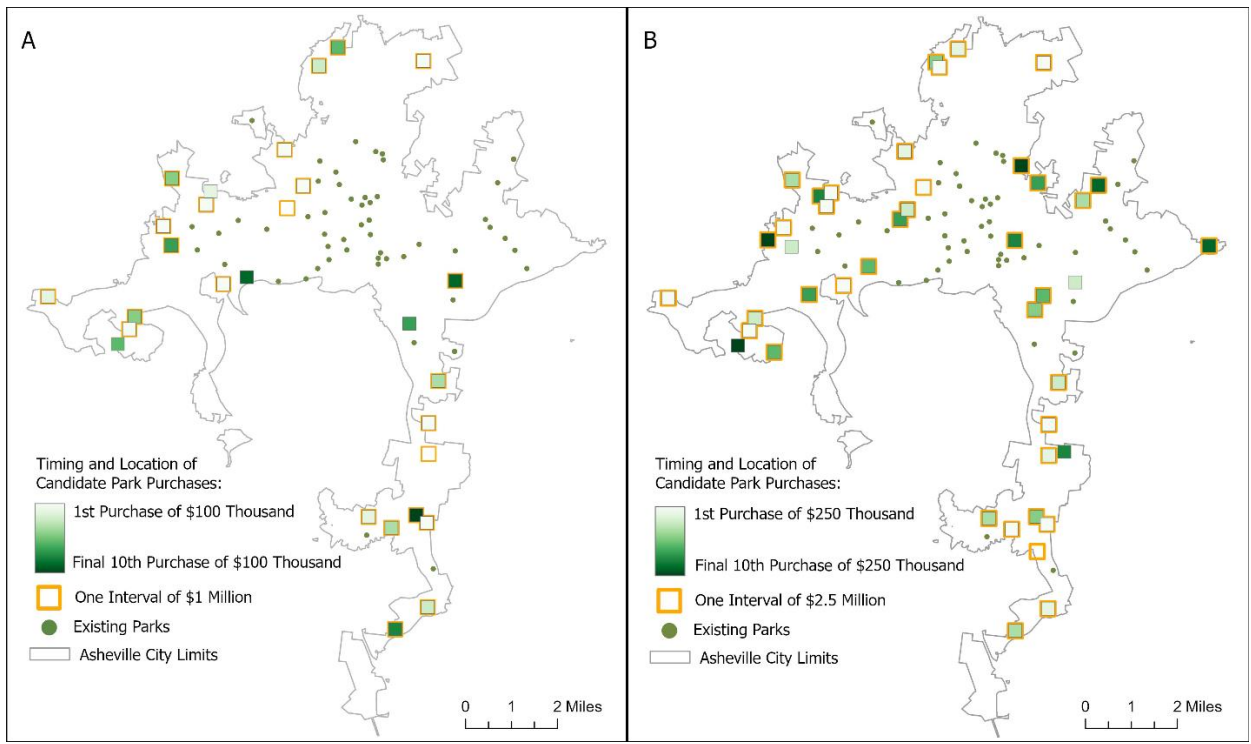

**Fig. 10.** Parks selected from long-term vs. myopic planning with $1 million budget (panel A) and $2.5 million budget (panel B)

parks in the most underserved areas of Asheville. In the initial years of myopic spending, the model locates new parks to produce a fairly even distribution of parks throughout the northern, western, and southern periphery. In later years of spending, the model revisits these regions to select additional candidate parks that further improve access; the model selects candidate parks that are physically nearer to the central and eastern regions. A total of 26 parks are selected in the myopic spending method. In contrast, new 24 candidate parks are selected with long-term planning. Among these, 22 parks are the same between the planning methods. Two parks are unique to the long-term planning process, and there are four unique parks with myopic planning.

With a higher total budget of $2.5 million, there is a similar pattern (Figure 10B). The number of uniquely selected candidate parks for long-term (one) vs. myopic planning (three) is lower when the budget is higher. We observe the result of this in the objectives as well; the access difference between planning strategies (long-term vs. myopic) is much narrower when the budget is higher (Figure 9). These results suggest that long-term planning is better, but a higher budget may be able to mitigate the effects of myopic planning.

**5.3 Strategic Demographic Emphases**

One approach to improve access is to target specific demographic groups. In this subsection, we analyze the effects of strategically emphasizing access for Black residents (Asheville Parks and Recreation, 2022). We use the *Min Max Dev Cap* model with a budget of $500,000. The comparative analyses are run with the strategic emphasis parameter set as $q_{Black} = 10$ (with strategic importance) vs. $q_{Black} = 1$ (without) with all others held to 1 ($q_r = 1, \forall r \in R \setminus \{Black\}$). The comparator of no strategic emphasis uses a value of 1 for each group $r \in R$. The value of 10 was calibrated as the threshold at which the optimal solution changed (tested weights from 0 to 50 in increments of 5).

Three new park locations differ based on whether the planners emphasize access for Black residents (Figure 11). With strategic emphasis on Black residents, the three unique parks are in the central and southern regions. Without strategic emphasis, the three unique parks are spread in the north, west, and southern regions. This occurs because of the distribution of Asheville residents. The northern region has a relatively large population across all race/ethnicities (Figure 3) but a relatively small number of Black residents.

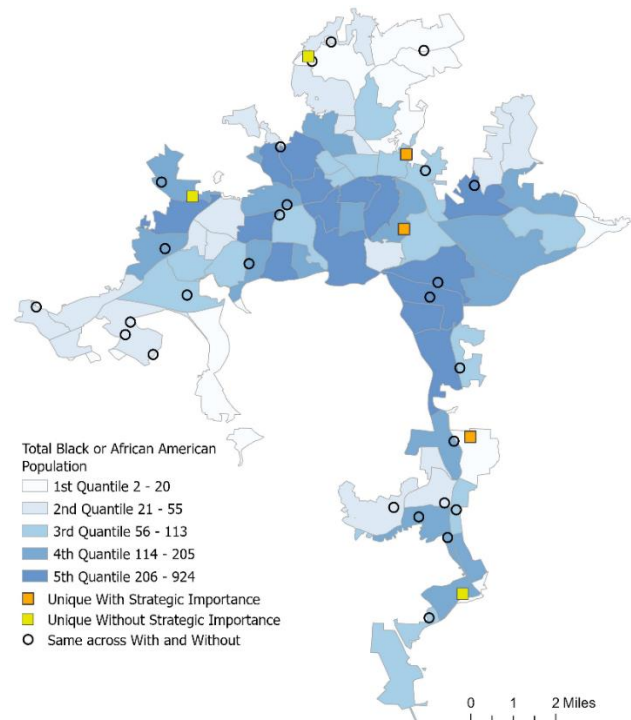

**Fig. 11.** Locations for new parks with vs. without strategic emphasis of Black residents

Most of the new parks are common to both analyses and are located throughout the city. In either analysis, new park locations are selected in some of the block groups with the largest populations of Black residents. One exception is in the central region. Few new parks are located in the center of the city because this area already has good park access (Figure 1). We note that with or without emphasis, access for Black residents improves vs. a baseline of $b = \$0$, though there are greater benefits when there is a particular emphasis.

**5.4 Access Thresholds**

Distance from residents to their primary parks is the highest priority category of access. While there is a target of 0.5 miles, in practice, this is a proxy for an individual's willingness to travel to a park. To analyze the stability of the locations across different thresholds, we evaluated alternate distances of 1 and 1.5 miles. We maintain a constant demographic weight of one for all racial demographics and a budget of $500,000. The majority of the primary parks remain consistent throughout all model instances, ("All Three Distances" in Figure 12). We note, however, that the primary parks added as a result of increased allowable distance from residents to parks tend to exist near the periphery of Asheville city limits.

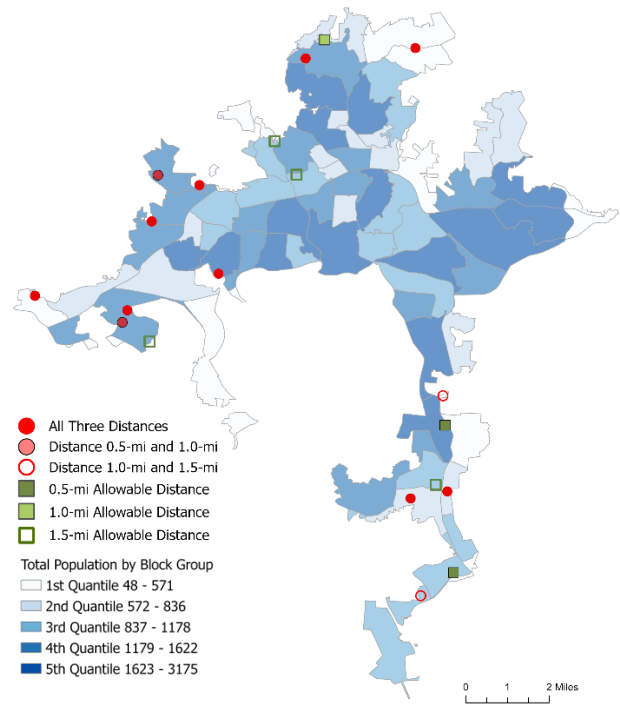

**Fig. 12.** Locations for new parks with vs. without strategic emphasis of Black residents

**5.5 Comparison with Alternative Location Optimization Models**

Finally, we compare access-oriented optimization with other approaches to location optimization. In Table 3, we present outcomes from the *Min Max Dev Cap* model (at budgets of $500,000 and $1 million) with three classic facility location formulations (p-median, p-center, and maximum coverage models). In the comparator models, the limit to the number of new parks (p) is set equal to the number of new parks in the optimal solution of the *Min Max Dev Cap* model and includes a side constraint of budget. Values presented in bold indicate the best value for each metric, stratified by budget category.

We make the following observations. At both budget levels, the access-oriented model (*Min Max Dev Cap*) outperforms the classic models in terms of weighted deviations (overall and worst-case group). It has moderate performance in demand-weighted distance, maximum distance, and coverage. The p-median models have the lowest demand-weighted distance (its objective function) with moderate-to-poor performance in weighted deviations. The p-center and maximum-coverage models lead in maximum distance and coverage, respectively, though not by a large margin; they perform very poorly in terms of access deviations. Note however that the allocation variables, $x_{kl}$, are not included in their objective functions, and these do not reflect optimized values. Overall, the *Min Max Dev Cap* models substantially improve overall access to parks and perform well according to metrics from other location models.

**Table 3.** Comparison of the Min Max Dev Cap model with p-median, p-center, and maximum covering models.

| Outcome | With $500,000 budget | | | | With $1 million budget | | | |
|---|---|---|---|---|---|---|---|---|
| | Min Max Dev Cap | p-Median (p=16) | p-Center (p=16) | Max Coverage | Min Max Dev Cap | p-Median (p=24) | p-Center (p=24) | Max Coverage |
| Overall weighted deviations, $\sum_{r \in R} \alpha_r$ | **463,203** | 580,775 | 921,690 | 1,781,450 | **407,452** | 502,507 | 938,120 | 2,262,275 |
| Weighted deviations for group with most, $\alpha^{max}$ | **349,687** | 442,269 | 697,007 | 1,388,542 | **307,844** | 382,595 | 710,120 | 1,742,181 |
| Demand-weighted distance (person-miles), $\sum_{l \in L} \sum_{r \in R} \sum_{k \in K} t_{lr} d_{kl} x_{kl}$ | 74,367 | **55,438** | 158,319 | 367,852 | 65,407 | **48,367** | 153,886 | 461,699 |
| Maximum distance from any location $l \in L$ to a selected park (miles) | 3.3 | **2.8** | **2.8** | 14.9 | 3.3 | 3.1 | **2.8** | 14.9 |
| Count of locations $l \in L$ that have any park located within 0.5 miles | 33 | 36 | 31 | **39** | 37 | 39 | 35 | **44** |
| Number of new parks selected | 16 | 13 | 12 | 18 | 24 | 21 | 20 | 23 |
| Number of new parks selected that are the same as in Min Max Dev Cap | -- | 5 | 7 | 10 | -- | 14 | 9 | 14 |

## 6. Discussion

Access to parks is fundamental to thriving communities, but it is well-established that there are disparities among demographic subpopulations, including by race and ethnicity (Rigolon, 2016). In this work, we provide a park planning optimization model to improve equitable access. The approach provides a framework to support planning efforts, and it goes beyond describing inequities to supporting steps to address them. In this section, we discuss the modeling and data collection (Section 6.1), policy analyses (Section 6.2), considerations for use (Section 6.3), and limitations (Section 6.4).

### 6.1 Modeling and Data Collection

The model considers four aspects of park access: distance, capacity, heat, and tree cover. Each has a range of values that are considered "good access," and deviations measure how far each resident is from good access by category. These deviations are weighted and aggregated together into a single metric: the sum of the weighted deviations by category for each resident in each location.

The framing of this metric of access is flexible. Other elements could be incorporated relatively easily by adding in their weighted deviations, and any of the current categories could be removed if not of interest. Corresponding constraints to define the deviations would be similarly revised. This flexibility is important because what constitutes good access is community-driven and not universal. In addition to the categories themselves, the ranges of allowable values can be adjusted to account for location variations. For example, in Tucson, Arizona, the allowable values for heat and tree cover would likely be very different than in mountainous, wooded Asheville. Alternative approaches could also consider multiple objectives. This could be framed as a bi-objective model that trades off worst-case ($\alpha^{max}$) with overall ($\sum_{r \in R} \alpha_r$) access. This trade-off would be particularly useful in settings where the objectives are in conflict; in Asheville, they produce similar solutions. It could also be structured as a multi-objective model in which each element of access is separately considered in its own objective.

The model captures deviations from good access in two ways. The first is to use pre-calculated deviations when they are exogenously known. For example, heat and tree cover deviations for each park $k \in K$ ($c_k^+, c_k^-, v_k^+, v_k^-$) do not change based on the decisions and can be pre-calculated outside of the model. The second is to use slack variables to measure the deviations that are endogenous to park location and primary park decisions. The distance deviations ($\tilde{d}_l^+$) for each resident location $l \in L$, and the number of people over capacity ($\tilde{a}_k^+$) at park $k \in K$ depend on the primary park designations. Note that capacity deviations are more difficult to model because they are endogenously counted (deviation variable $\tilde{a}_k^+$) and then penalized for each resident associated with the park (decision variable $x_{kl}$), i.e., the nonlinear term $\sum_{l \in L} \sum_{k \in K} t_{lr} \tilde{a}_k^+ x_{kl}$ in constraints (2). This is a feature of *person*-centered access because the deviation is experienced by every resident that is associated with the park. In contrast, a conventional *location*-centered approach would only penalize $\tilde{a}_k^+$ once for each park, e.g., $\sum_{k \in K} \tilde{a}_k^+$. We suggest that while modeling person-centered experience is more complex (i.e., must be linearized), it is more reflective of practice, where each person's individual experience of parks affects their access.

Extensive geospatial datasets were compiled to support these analyses (Appendix Table B1). The candidate locations were derived from city parcels using zoning data and maps of existing infrastructure. This set could be further enhanced through feedback from local residents and projections of property availability. Details on park access metrics were derived from geospatial data, including networks to estimate walking distances. Further work could incorporate public transportation and more individualized thresholds of access. Population sizes by race/ethnicity were determined from the 2020 Census (*Race (P1)*, 2020). A recent report noted potential over- and under-counting (US Census Bureau, 2022), and future work could consider uncertainty in population sizes.

Analysis of the model structure indicates that there may be limits to the ability of a planner to use new park investments to improve access (Proposition 1) and that adding parks (Proposition 2) or budget (Corollary) will not reduce access. In future work, these results could be incorporated into more advanced solution approaches, e.g., via submodular optimization (Küçükyavuz & Yu, 2023).

**6.2 Policy Analyses**

There are multiple mechanisms to address inequities, and we use the model to evaluate several policies. One option is the budget itself. We find, intuitively, that increasing the budget decreases weighted deviations from good access, i.e., increases access. This trend occurs with either type of objective function – increasing the budget improves access either when the goal is to improve access for the group that has the worst access (objective function 1) or to improve for all residents (objective function 20). We observe decreasing marginal benefits, however. For each incremental increase in available budget, the amount by which access improves decreases. Eventually weighted deviations hit a threshold and cannot be improved by location decisions, as suggested by Proposition 1. Because there are often limited budgets for the acquisition of new parks, these results suggest that park planners consider these trade-offs in the broader context of the portfolio of spending to improve access. Other options to improve access could include improved programming, new amenities, and/or renovations at existing parks. In general, though, we find that adding a park anywhere in the city will not hurt access and may improve it (per Proposition 2). This suggests that if the city obtains access to additional land in an underserved area, e.g., through a new easement program or a bequeathment, it may be worthwhile for park access to develop a park on it.

For a given budget, one way to improve access is to reconsider the planning process. In a myopic approach, new park(s) may be selected each year; for a 10-year horizon, this means plans are made 10 times. In a long-term approach, new parks are selected all at once; there is a single plan. In both myopic and long-term term planning, the same amount of money is available; it is the timing of the decisions that change. We find that, while both planning horizons can substantially improve access vs. baseline, a long-term approach can amplify the benefits of the available budget. A long-term plan may improve access by 5% with a budget of $1 million and 1% if the total is $2.5 million compared to a myopic plan. We do not suggest that 10 years is a magic number but rather that there are benefits to *coordinated decisions*. This insight has implications for the development of other software to reduce park inequities (e.g., ParkServe, The Trust for Pulic Land, 2022). These programs designate areas within cities that have low access, and they do not generally allow users to evaluate how access changes if multiple parks are added concurrently.

Cities may also choose to prioritize particular subpopulations when addressing inequities. We study the effects of strategically emphasizing access for Black residents of Asheville. In this case, we observe that new parks are located in areas close to where Black residents live, and people in regions with larger Black populations do not have to travel as far to visit their primary parks. For residents in areas with smaller numbers of Black residents, including the Black residents of these areas themselves, the distance to the primary parks is greater than if there were no emphasis. We note, though, that with or without strategic emphasis on a particular subpopulation, when money is allocated to purchase new parks, access improves both for the subpopulation and overall as compared to baseline.

We also consider how the distance threshold affects park locations. Throughout the other analyses, we use 0.5 miles between a resident location and park as the ideal distance threshold. When this is increased to 1 or 1.5 miles, many of the parks are located in the same locations. Among the few instances that parks are selected in different locations than at the baseline of 0.5 miles, they are located in more isolated regions. This suggests that the results are at least to some degree robust across different thresholds of access. This is important because access looks different for everyone. A 20-minute walk may be reasonable for some but not for others. Park locations need to consider the inclusion of all residents, and while single thresholds do not capture this variability, the stability in locations selected suggests they may be a decent proxy when more detailed data is unavailable.

**6.3 Considerations for Model Use**

A goal of this work is to support local community efforts to increase access to parks. The model was developed with feedback from Asheville park planners, and we met with APRD during the development stages and presented results after the analyses were conducted. However, models are inherently simplifications of reality, and carefully considering locations in partnership with community members is crucial.

In particular, there are a handful of features that warrant further consideration. One is related to gentrification. Parks may be built to serve particular communities, but new parks can be a draw for wealthier residents to move to the area (Eldridge et al., 2019). Property values may increase, and the original community that the park was intended to serve may no longer be able to afford to live there. In Chicago, the 606 Trail was installed as a multi-use recreational park but has led to extensive displacement of long-time residents (Black, 2020). Similarly, the opening of the High Line walkway in New York led to rent increases of 68% in the neighboring areas between 2009-2013 and 2014-2018 (Black and Richards, 2020). In these cases, not only did the new park not help the original community, the new park harmed it. To reduce the risks of gentrification induced by new parks, park and urban planners could partner with housing organizations to support affordable housing in these areas (Eldridge et al., 2019; Rigolon and Németh, 2018).

We also note that adding parks is necessary but not sufficient to address inequities. Investments in existing parks, including renovations to ensure accessibility and to add resident-desired amenities, may further improve well-being (Ayala-Azcárraga et al., 2019). We note that there is a growing literature on greening by which communities have ownership and bring their cultures into the art and design of the park space (Jeong et al., 2021). As with much applied operations research, the work itself does not end with a model; rather, this work is one step in the process toward improved access.

This modeling approach may be able to be applied in other cities that seek to locate new parks to improve equity. The data itself is likely to be available, e.g., in the US via publicly available datasets and from city partners who often have a GIS contact on staff. The specific elements of access could be tailored to the local community (e.g., to include other features such as accessibility or safety, Ayala-Azcárraga et al., 2019; Williams et al., 2020)). Long-term plans that have annual budgets could be optimized by stratifying park selection and budgets by year. In general, the structure may be versatile. The primary considerations on whether it would generalize to a particular city would be whether it has similar goals to those in Asheville, i.e., is the city seeking to invest in new parks and do its park planners aim to locate new parks through the lens of equity? We note that there is a growing recognition of the importance of addressing the status quo of inequitable access (National Recreation and Park Association, 2011; *Padilla, Collins Introduce Bipartisan Outdoors for All Act to Fund Urban Parks Projects*, 2023). We expect interest in these questions to continue to grow. Beyond the focus of this paper, the model may also be able to be adapted to a park closure setting, e.g., to identify the park that does the least harm to equitable distribution of resources.

### 6.4 Limitations

There are limitations with regard to both our data and models. Concerning data, we were unable to acquire information defining population density *within* blocks groups. Therefore, we assume homogeneity of population distribution within each resident location $l \in L$. The lack of population density information decreases the accuracy of the representation of demographic data. In the calculation of distance from residents to parks, we set the origin of each resident location as the block group geometric center-point rather than the center of population density, which was unavailable. Due to the lack of visitation data, the capacitated models take the perspective of a municipal planner and do not include resident choice. To mitigate the effects we compare the capacitated model throughout with an uncapacitated model that is not affected by this assumption.

### 7. Conclusion

This paper focuses on the development of a new integer programming model to support efforts to improve access to urban recreation. We provide a park location tool that integrates core demographic, geospatial, fiscal, and environmental factors that affect person-centered park access. Using Asheville, NC as a case study, we complete extensive data collection to translate current-state park realities into usable model inputs and conduct analyses to answer key policy questions of budget use, strategic priorities, and metrics of access. The use of publicly available data and reasonable solve times suggest the potential to scale to other cities that seek to locate new parks to improve equity. This research is one step towards improving access in recreation and provides a range of open questions.

Future work could consider the effects of uncertainty in long-term planning, resident-choice models, as well as other demographic stratifications such as socioeconomic status, disability, gender, and age. The latter could lead to further study on park accessibility (e.g., wheelchair accommodations) and individualized measures of access, including the difficulty of terrain while walking to parks (e.g., hills). Resident-choice models could use gravity-based approaches (Suárez-Vega et al., 2012; Y. Zhang & Atkins, 2019) based on distance or park features (Costigan et al., 2017; Kaczynski et al., 2014) to estimate visitation. The growing availability of data from social media (Chuang et al., 2022; Hamstead et al., 2018; Tenkanen et al., 2017) and cell phones (He et al., 2022; Zhai et al., 2018) also has promise to estimate visitation to parameterize resident-choice models. Taken together, the presented model and the many future directions highlight the potential for operations research to support decisions in recreation.

Rec-Underserved-Areas.pdf

Neema, M. N., & Ohgai, A. (2010). Multi-objective location modeling of urban parks and open spaces: Continuous optimization. *Computers, Environment and Urban Systems*, *34*(5), 359–376. https://doi.org/10.1016/j.compenvurbsys.2010.03.001

Niedermann, B., Oehrlein, J., Lautenbach, S., & Haunert, J. H. (2018). A network flow model for the analysis of green spaces in urban areas. *Leibniz International Proceedings in Informatics, LIPIcs*, *114*(13), 1–13. https://doi.org/10.4230/LIPIcs.GIScience.2018.13

Öhman, K., & Lämås, T. (2005). Reducing forest fragmentation in long-term forest planning by using the shape index. *Forest Ecology and Management*, *212*(1–3), 346–357. https://doi.org/10.1016/j.foreco.2005.03.059

Owen, S. H., & Daskin, M. S. (1998). Strategic facility location: A review. *European Journal of Operational Research*, *111*(3), 423–447. https://doi.org/10.1016/S0377-2217(98)00186-6

*Padilla, Collins Introduce Bipartisan Outdoors for All Act to Fund Urban Parks Projects*. (2023). Press Release. https://www.padilla.senate.gov/newsroom/press-releases/padilla-collins-re-introduce-bipartisan-outdoors-for-all-act-to-fund-urban-parks-projects

Pourrezaie-Khaligh, P., Bozorgi-Amiri, A., Yousefi-Babadi, A., & Moon, I. (2022). Fix-and-optimize approach for a healthcare facility location/network design problem considering equity and accessibility: A case study. *Applied Mathematical Modelling*, *102*. https://doi.org/10.1016/j.apm.2021.09.022

*Race (P1)*. (2020). United States Census Bureau.

Rigolon, A. (2016). A complex landscape of inequity in access to urban parks: A literature review. *Landscape and Urban Planning*, *153*, 160–169. https://doi.org/10.1016/j.landurbplan.2016.05.017

Rigolon, A. (2017). Parks and young people: An environmental justice study of park proximity, acreage, and quality in Denver, Colorado. *Landscape and Urban Planning*, *165*, 73–83. https://doi.org/10.1016/j.landurbplan.2017.05.007

Rigolon, A., Fernandez, M., Harris, B., & Stewart, W. (2019). An Ecological Model of Environmental Justice for Recreation. *Leisure Sciences*. https://doi.org/10.1080/01490400.2019.1655686

Rigolon, A., & Németh, J. (2018). "We're not in the business of housing:" Environmental gentrification and the nonprofitization of green infrastructure projects. *Cities*, *81*. https://doi.org/10.1016/j.cities.2018.03.016

Sefair, J. A., Molano, A., Medaglia, A. L., & Sarmiento, O. L. (2012). Locating Neighborhood Parks with a Lexicographic Multiobjective Optimization Method. In M. Johnson (Ed.), *Community-Based Operations Research*. Springer Science+Business Media. https://doi.org/0.1007/978-1-4614-0806-2_6

**Supplementary Materials**

"An Optimization Approach to Improve Equitable Access to Local Parks"

In this collection of appendices, we provide the following supplementary information. Appendix A presents linearization of the constraints. Appendix B provides additional information related to data collection and preparation of model inputs. Appendix C presents supplemental results.

**Appendix A**

Constraints (2), (9), and (11) are nonlinear because they include the product of decision variables. We define additional parameters and decision variables and present the linearized constraints as follows.

| | |
|---|---|
| ***Parameters*** | |
| $\mu^{cap+}$ | big M (max) value for over capacity of parks |
| $\mu^{maxdist}$ | big M (max) value for distance from resident locations to parks |
| $\mu^{maxcap}$ | big M (max) value for capacity of parks |
| ***Decision Variables*** | |
| $\pi_{kl}^{cap+}$ | linearization variable for the over capacity of park $k \in K$ experienced by location $l \in L$ |
| $\pi_l^{actdist}$ | linearization variable for distance slack for location $l \in L$ |
| $\pi_k^{actcap}$ | linearization variable for limiting capacity slack for park $k \in K$ |

*Linearization of Park Access Deviations, Constraints (2):*

To linearize Constraints (2), we define $\pi_{kl}^{cap+} \coloneqq \tilde{a}_k^+ x_{kl}$ and replace (2) with an exact reformulation: (22-26):

$$\alpha_r = \sum_{l \in L} \sum_{k \in K} \left( q_r t_{lr} \left[ \begin{pmatrix} n^{dist} w^{dist+} \tilde{d}_l^+ \\ +n^{cap} w^{cap+} \pi_{kl}^{cap+} \end{pmatrix} + \begin{pmatrix} n^{heat} w^{heat+} c_k^+ \\ +n^{heat} w^{heat-} c_k^- \\ +n^{tree} w^{tree+} v_k^+ \\ +n^{tree} w^{tree-} v_k^- \end{pmatrix} x_{kl} \right] \right) \quad \forall r \in R \quad (22)$$

$$\pi_{kl}^{cap+} \leq \mu^{cap+} x_{kl} \quad \forall k \in K, l \in L \quad (23)$$

$$\pi_{kl}^{cap+} \leq \tilde{a}_k^+ \quad \forall k \in K, l \in L \quad (24)$$

$$\pi_{kl}^{cap+} \geq \tilde{a}_k^+ - (1 - x_{kl}) \mu^{cap+} \quad \forall k \in K, l \in L \quad (25)$$

$$\pi_{kl}^{cap+} \geq 0 \quad \forall k \in K, l \in L \quad (26)$$

*Linearization of Distance Slack Calculation, Constraints (9):*

To linearize Constraints (9), we define $\pi_l^{actdist}$ and replace (9) with an exact reformulation: (27-31).

$$\tilde{d}_l^+ - \sum_{k \in K} d_{kl} x_{kl} + m + \pi_l^{actdist} - u_l^{dist} m \leq 0 \quad \forall l \in L \tag{27}$$

$$\pi_l^{actdist} \leq \mu^{maxdist} u_l^{dist} \quad \forall l \in L \tag{28}$$

$$\pi_l^{actdist} \leq \sum_{k \in K} d_{kl} x_{kl} \quad \forall l \in L \tag{29}$$

$$\pi_l^{actdist} \geq \sum_{k \in K} d_{kl} x_{kl} - (1 - u_l^{dist}) \mu^{maxdist} \quad \forall l \in L \tag{30}$$

$$\pi_l^{actdist} \geq 0 \quad \forall l \in L \tag{31}$$

*Linearization of Capacity Slack Calculation (Constraint 11):*

To linearize Constraints (11), we define $\pi_l^{actcap}$ and replace (11) with an exact reformulation: (32-36)

$$\tilde{a}_k^+ - \sum_{l \in L} \sum_{r \in R} t_{lr} x_{kl} + a_k + \pi_k^{actcap} - u_k^{cap} a_k \leq 0 \quad \forall k \in K \tag{32}$$

$$\pi_k^{actcap} \leq \mu^{maxcap} u_k^{cap} \quad \forall k \in K \tag{33}$$

$$\pi_k^{actcap} \leq \sum_{l \in L} \sum_{r \in R} t_{lr} x_{kl} \quad \forall k \in K \tag{34}$$

$$\pi_k^{actcap} \geq \sum_{l \in L} \sum_{r \in R} t_{lr} x_{kl} - (1 - u_k^{cap}) \mu^{maxcap} \quad \forall k \in K \tag{35}$$

$$\pi_k^{actcap} \geq 0 \quad \forall k \in K \tag{36}$$

# Appendix B

We present a table of data sources (Table B1), additional maps (Figures B1-B4), and a discussion of alternative approaches to prepare the data.

**Table B1.** Data collected and data sources for model inputs.

| Model Element | Data | Data Description | GIS Type | Year | Source | Reference |
|---|---|---|---|---|---|---|
| General boundaries | Asheville city limits | Spatial representation of Asheville city limits | polygon | 2017 | City of Asheville | (*Asheville City Limits*, 2017) |
| Set of resident locations, $L$ | Census Block Groups | Spatial mapping and non-spatial list of Block Groups | polygon | 2019 | US Census | (US Census Bureau, 2021) |
| Set of races/ethnicities, $R$ | Racial/ethnic classifications | Non-spatial list of racial/ethnic classifications | n/a (tabular) | 2020 | US Census | (*Race (P1)*, 2020) |
| Population counts, $t_{lr}$ | Race/ethnicity | Number of individuals in each racial/ethnic group (2020 Block Groups) | n/a (tabular) | 2020 | US Census | (*Race (P1)*, 2020) |
| | Census Block Groups | Spatial distribution of Census Block Groups | polygon | 2020 | US Census | (Block Groups 2020) |
| Set of existing parks, $K^{existing}$ | Existing parks | Spatial distribution of existing parks | polygon | 2021 | City of Asheville | (*Asheville Parks*, 2019) |
| Set of candidate park locations, $K^{candidate}$ | Buncombe parcels | Spatial distribution of land parcels | polygon | 2020 | Buncombe County | (Buncombe County GIS, 2023) |
| | Building zoning codes | Spatial distribution of zoning districts | polygon | 2020 | City of Asheville | (*Asheville Zoning*, 2020) |
| | Floodways | Spatial distribution of flood zones | polygon | 2021 | Federal Emergency Management Agency (FEMA) | (*FEMA Flood Map*, 2021) |
| | Water features | Spatial distribution of water – lakes, ponds, streams | polygon | 2019 | US Census | (*Water*, 2019) |
| Park cost, $f_k$ | Land values | Land value for each park | n/a (tabular) | 2020 | City of Asheville, Buncombe County | (*Asheville Parks*, 2019; Buncombe County GIS, 2023) |
| Park Capacity, $a_k$ | Residents per acre | Standard number of residents per park acre | n/a (tabular) | 1965 | American Planning Association | (Moeller, 1965) |
| | Park size | Size of park in acres | polygon | 2020 | Environmental Systems Research Institute (ESRI) | (ESRI, 2023) |
| Park Heat, $c_k^+, c_k^-$ | Heat severity | Index (1-5) of the severity of heat above the city average within a cell | raster | 2021 | Trust for Public Land | (The Trust for Public Land, 2022) |
| Park Tree cover, $v_k^+, v_k^-$ | Tree cover | Percentage (0-100) of tree cover within a cell | raster | 2016 | Multi-Resolution Land Characteristics (MRLC) Consortium | (Multi-Resolution Land Characteristics Consortium, 2020) |
| Walking distance, $d_{kl}$ | Streets | Network of Asheville streets that allow pedestrian travel | line | 2020 | ESRI | (ESRI, 2023) |

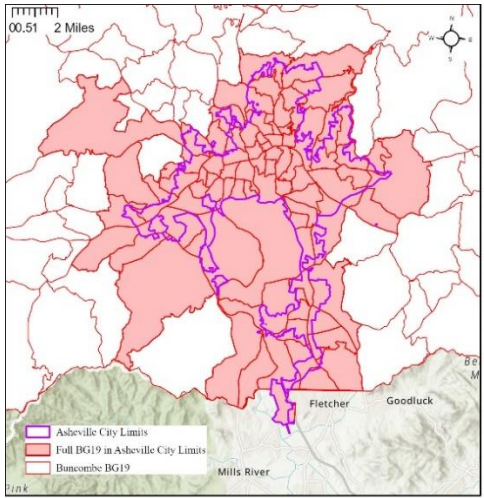

**Fig. B1.** Block groups (2019) and Asheville City Limits.

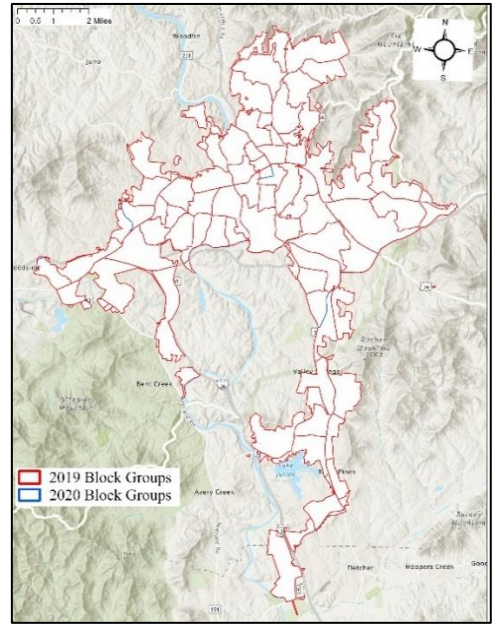

**Fig. B2.** Redistricting of Block Groups from 2019 to 2020

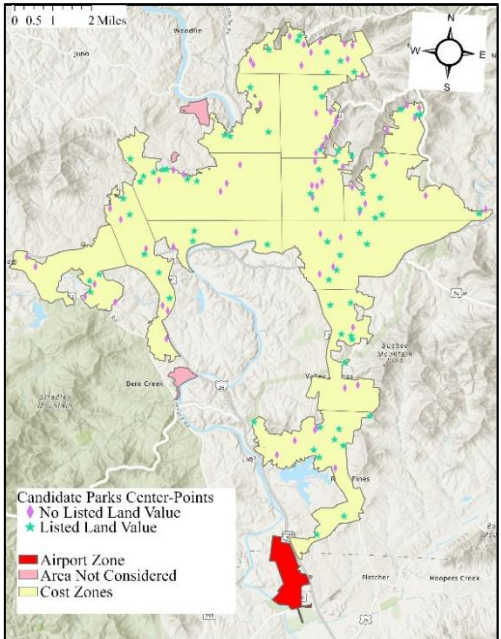

**Fig. B3.** Land values and cost zones

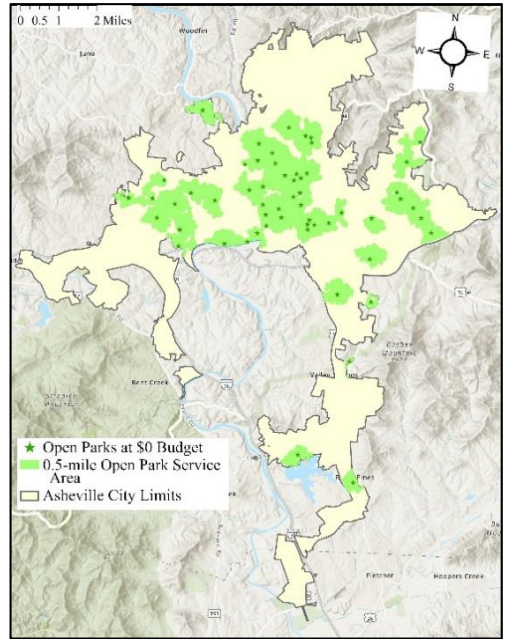

**Fig B4.** Budget of $0, park locations

**Discussion of Application of Geospatial Methods**

In other settings, it may be useful for researchers and practitioners to select alternative methods for the use of the geospatial data. If edge effects are prominent, e.g., there may be parks outside of city limits that residents may consider primary, the modeler may choose to expand set $K^{existing}$ to include these parks. If park planners aim to consider all residents within block groups, including those beyond city limits, the modeler should not clip the block groups to exclude areas outside of city limits (Figure B1).

In this paper, we calculated the distances $d_{kl}$ between the centroid of resident location $l \in L$ and *centroid* of candidate park $k \in K$ (Section 4.8). Modelers could alternatively parameterize $d_{kl}$ to represent the distance between the centroid of each location $l \in L$ and the nearest *entrance* to park $k \in K$. If a park has multiple entrances and/or open boundaries, we suggest using the minimum distance for each resident location $l \in L$.

More granular population datasets may also be available depending on analysis goals. The US Census data has been updated in 2023, after the initial analysis, and is now available for 2020 at the more granular block level (U.S. Census Bureau, 2020). Population estimates at the 100m resolution, unstratified by race/ethnicity, are available from WorldPop (Bondarenko et al., 2020).

**Appendix C**

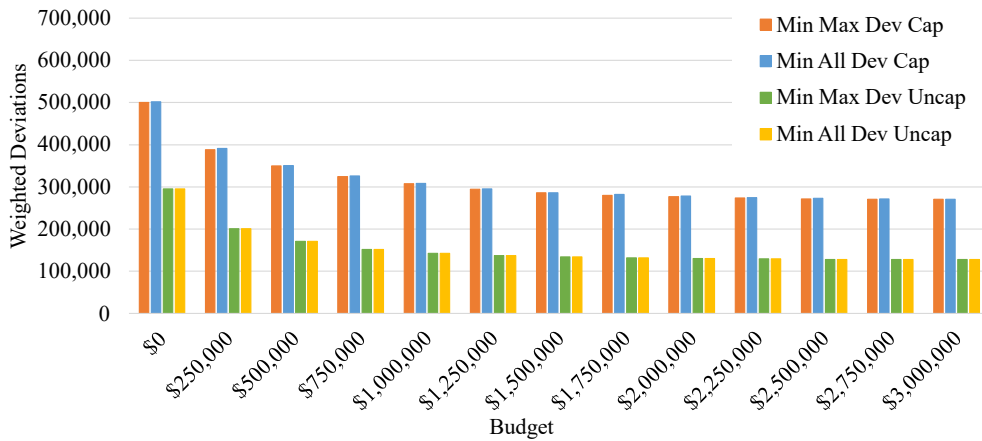

**Fig. C1.** Weighted deviations for subpopulation with the highest

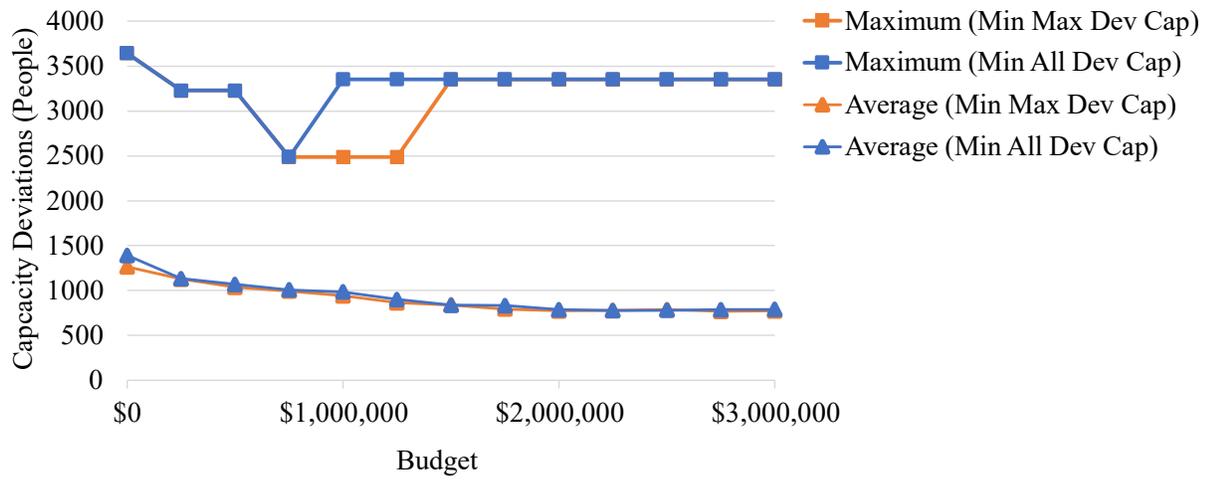

**Fig. C2.** Resident-experienced capacity deviations

**References for Supplementary Materials**